\documentclass[a4paper, twosided, 10pt,reqno]{amsart}
\usepackage[hidelinks]{hyperref}
\usepackage{amsthm,amsmath,amssymb,mathtools}
\usepackage{enumitem}
\usepackage[english]{babel}

\usepackage{enumitem}
\usepackage{microtype}
\usepackage{graphicx}
\usepackage{multicol}
\usepackage{tikz-cd}
\usetikzlibrary{patterns}
\usepackage{pgfplots}
\usetikzlibrary{intersections}
\usepackage[dvipsnames]{xcolor}

\newcommand*{\defeq}{\mathrel{\vcenter{\baselineskip0.5ex \lineskiplimit0pt
			\hbox{\scriptsize.}\hbox{\scriptsize.}}}%
	=}
	
\newcommand*{\eqdef}{=%
	\mathrel{\vcenter{\baselineskip0.5ex \lineskiplimit0pt
			\hbox{\scriptsize.}\hbox{\scriptsize.}}}}
		
\renewcommand{\epsilon}{\varepsilon}
\renewcommand{\theta}{\vartheta}
\renewcommand{\phi}{\varphi}
\newcommand{\R}{\mathbb{R}}

\newcommand{\N}{\mathbb{N}}

\newcommand{\esse}{\mathbb{S}}

\usepackage{accents}

\makeatletter

\newcommand{\tr}{\text{tr}}

\title{The exterior Dirichlet problem for the Weingarten equation}
\author{Angelo Benedetti}

\address{Universit\`a degli Studi dell'Aquila \\
		Dipartimento di Ingegneria e Scienze dell'Informazione e Matematica\\
		via Vetoio 1 \\
		67100 L'Aquila\\
		Italy}
\email{angelo.benedetti@graduate.univaq.it}

\newtheorem{lemma}{Lemma}[section]
\newtheorem{theorem}[lemma]{Theorem}
\newtheorem{theoremletter}{Theorem}

\newtheorem{lemmaletter}[theoremletter]{Lemma}

\newtheorem{corollary}[lemma]{Corollary}
\newtheorem{proposition}[lemma]{Proposition}
\theoremstyle{definition}
\newtheorem{definition}[lemma]{Definition}
\newtheorem*{definition*}{Definition}
\theoremstyle{remark}
\newtheorem{remark}[lemma]{Remark}
\newtheorem*{remark*}{Remark}

\begin{document}
	\maketitle
	\begin{abstract}
	We solve the exterior Dirichlet problem for the Weingarten equation of minimal type over the complementary of domains satisfying an interior sphere condition, and investigate some related problems. Namely, by linearizing the Weingarten equation, we give an asymptotic expansion for exterior solutions and we completely identify which classes of Weingarten minimal ends satisfy a maximum principle at infinity.
	Building on the previous results, we study symmetries of the solutions and questions related to the uniqueness of the solutions.
	\end{abstract}
	\section{Introduction}
Let $\Sigma\subset\R^3$ be an immersed oriented surface with principal curvatures \(k_1\), \(k_2\), with \(k_1\geq k_2\). We say that $\Sigma$ is a \textit{Weingarten surface} when it satisfies at every point the \textit{Weingarten equation}:
	\begin{equation}\label{equation-weingarten}
		k_2=f(k_1),
	\end{equation}
	with \(f\in C^{1,1}_{\text{loc}}[c,+\infty)\) for some \(c\geq0\) such that \(f(c)=c\). The Weingarten equation is called \textit{uniformly elliptic} when:
	\begin{equation}\label{equation-uniformellipticity}
		-\frac{1}{\Lambda}\leq f'(t)\leq-\Lambda\quad\forall t\geq c,
	\end{equation}
	for some constant $\Lambda\in(0,1)$.
	The behavior of an elliptic Weingarten surface depends strongly on the value of \(c\). We say that the Weingarten equation is of \textit{minimal type} when:
	\begin{equation}\label{equation-minimaltype}
		c=0.
	\end{equation}
	The definition comes from the fact that Weingarten surfaces of minimal type are seen to share many features of minimal surfaces. In fact, the theory of uniformly elliptic Weingarten surfaces is, in a sense, the fully nonlinear analogue of constant mean curvature (cmc) surfaces and minimal surfaces. It is then interesting to investigate to what extent does the analogy hold, and try to develop new tools to prove the analogues of the main global theorems of cmc surfaces theory in this nonlinear setting. That is because, usually, the classical tools of cmc surfaces theory and minimal surfaces theory do not longer hold: we do not have, for instance, a Weierstrass representation, nor an underlying variational problem. Thus, one has to find new paths to get to familiar results, which shed new light and bring more insight to the whole theory. For example, Rosenberg, Sa Earp \cite{saearpweingarten} and Aledo, Espinar, Gálvez \cite{espinar} extended part of the Korevaar-Kusner-Meeks-Solomon theory to Weingarten surfaces of cmc type; Nelli, Pipoli, Tassi \cite{barbaraweingarten} showed a Ros-Rosenberg type theorem for Weingarten surfaces of cmc type; Fernández, Gálvez, Mira \cite{galvezquasiconformal} proved an analogue of the Bernstein theorem for Weingarten multigraphs.
	
	Following this line, Barbieri, Gálvez, Lian, Zhang recently described the asymptotic behavior of the annular ends of a uniformly elliptic Weingarten surface of minimal type with finite total curvature \cite{galvezbarbieri2026asymptotic}. The behavior coincides in a sense with the behavior of the rotational examples, which have been classified in \cite{fernandezweingarten}, without assuming uniform ellipticity, and it is seen to strongly depend on the value of \(f'(0)\) (see Theorem \ref{theorem-galvez-c}). Indeed, uniformly elliptic Weingarten multigraphs of minimal type satisfy curvature estimates depending on the distance from the boundary (see \cite[Theorem 5.2]{galvezquasiconformal}, \cite[Theorem 7.3]{galvezbarbieri2026asymptotic}). Thus, close to infinity, we have that the principal curvatures go to zero and, heuristically, the behavior of the Weingarten surface should depend from the behavior of the graph of \(f\) close to the origin.
	
	Observe that the analogous description for minimal surfaces is usually proved using the Weierstrass representation which, as already observed, is not available in this setting. Therefore, the authors had to study the asymptotic behavior of solutions of the associated PDE over \textit{exterior domains}.
	
	In this paper, we further pursue this line of investigation by proving existence of solutions to the Dirichlet problem over exterior domains, which is a problem already suggested in \cite{galvezbarbieri2026asymptotic}, and some related results.
	
	We assume the graphs to be oriented by the upper normal, and study solutions to \eqref{equation-weingarten}-\eqref{equation-uniformellipticity}-\eqref{equation-minimaltype} with limit unit normal \(e_3\) (exterior solutions to the uniformly elliptic Weingarten equation have a limit unit normal, as explained in \cite{galvezbarbieri2026asymptotic}) on \(\overline{\Omega}^C\), where $\Omega\subset\R^2$ is a bounded domain such that \(\overline{\Omega}^C\) is connected. We first show that, under a suitable assumption on the geometry of $\Omega$, we can solve the Dirichlet problem on $\overline{\Omega}^C$ with zero boundary datum (Section \ref{section-exteriordirichlet}). The proof relies on the Perron method, which requires some preparation in order to be applied to the fully nonlinear setting. This is done in Section \ref{section-dirichlet}, where we gather some known results on \(C^{2,\alpha}\) estimates for the uniformly elliptic Weingarten equation. Then, in Section \ref{section-perron}, we describe the Perron method for the uniformly elliptic Weingarten equation and the geometric barriers we will use.
	
	Section \ref{section-maximumprincipleinfinity} concerns the asymptotic expansion for uniformly elliptic Weingarten ends of minimal type. We remark that asymptotic expansion was first proved in \cite{galvezbarbieri2026asymptotic} for \(f'(0)=-1\). In this case, the underlying PDE is of class \(C^{1,1}\), and the proof goes by linearizing the equation. Nevertheless, for \(f'(0)\neq-1\) the Weingarten equation is not differentiable at umbilical points (i.e. points where \(k_1=k_2\)), and so we must prove that, close to infinity, every point of the graph of an exterior solution is sufficiently far from being umbilical (Lemma \ref{lemma-differentiability}). From the asymptotic expansion, various results follow. First of all, we completely identify which classes of Weingarten equations admit a maximum principle at infinity (Corollary \ref{corollary-maximumprincipleatinfinity} and Remark \ref{remark-maximumprincipleatinfinity}). Then, we give a comparison criterion for ends in the same Weingarten class which hold for every class of uniformly elliptic Weingarten ends of minimal type (Theorem \ref{theorem-weakmaximumprincipleatinfinity}). Using this result, we show that uniformly elliptic Weingarten ends of minimal type inherit the symmetries of their boundary data.
	
	Section \ref{section-krust} contains an extension to the Weingarten setting of a result of Krust for the minimal surface equation \cite{krustplateau}. We show that, when two distinct solutions \(u_a\) and \(u_b\) to the same Dirichlet problem exist, then there is a one-parameter family of solutions to the same problem foliating the region bounded by the graphs of \(u_a\) and\,\,\(u_b\).
	\section{Interior estimates for the Dirichlet problem}\label{section-dirichlet}
	Let $\Omega$ be an open domain of \(\R^2\). Given a continuous function \(\phi:\partial\Omega\to\R\), we remember that we refer as \(\mathcal{P}(\Omega,\phi)\) to the problem of finding a \(C^2\) function solving \eqref{equation-weingarten}-\eqref{equation-uniformellipticity}-\eqref{equation-minimaltype} in\,\,$\Omega$ and having \(\phi\) as boundary datum. For a strictly convex domain $\Omega$, the interior Dirichlet problem $\mathcal{P}(\Omega,\phi)$ was first solved in \cite{galvezbarbieri2026asymptotic}, under regularity assumptions on $\phi$:
	\begin{theoremletter}[\cite{galvezbarbieri2026asymptotic}, Theorem 1.1]\label{theorem-galvezdirichlet}
		Let $\Omega \subset \mathbb{R}^2$ be a strictly convex, bounded domain of class
		$C^{2,\alpha}$, and let $\varphi \in C^{2,\alpha}(\partial\Omega)$, with
		$\alpha\in(0,1)$. Then there exists a unique solution
		$u \in C^{2,\beta}(\overline{\Omega})$ (for some $0 < \beta \leq \alpha$)
		to $\mathcal{P}(\Omega,\phi)$, where
		$f \in C_{\mathrm{loc}}^{0,1}[0,+\infty)$.
		Moreover,
		\begin{equation}
			\label{eq:1.4}
			\lVert u \rVert_{C^{2,\beta}(\overline{\Omega})} \leq C,
		\end{equation}
		where $C$ depends only on
		\[
		\Lambda,\ \alpha,\ \max_{\partial\Omega}\kappa_{\partial\Omega},\
		\min_{\partial\Omega}\kappa_{\partial\Omega},\
		\lvert\partial\Omega\rvert_{C^{2,\alpha}},\
		\text{and }\lVert\varphi\rVert_{C^{2,\alpha}(\partial\Omega)}.
		\]
		Here, $\kappa_{\partial\Omega}$ denotes the curvature of $\partial\Omega$.
	\end{theoremletter}
Observe that the theorem applies in principle only for \(C^{2,\alpha}\) boundary datum. Nevertheless, in order to use the Perron method, we will need to be able to solve the problem for just continuous boundary datum. This is what we aim to prove in the rest of this section (Theorem \ref{theorem-dirichletC0}). The idea is the following.

Given a continuous function $\phi$ defined on $\partial\Omega$, with \(\partial \Omega\) of class \(C^{2,\alpha}\), we can approximate $\phi$ uniformly using a sequence $\{\phi_n\}_n\subset C^{2,\alpha}(\partial\Omega)$:
\[\phi_n \xrightarrow{C^0} \phi.\]
Then, for each \(n\in\N\), the problem $\mathcal{P}(\Omega,\phi_n)$ admits a unique solution \(u_n\in C^{2,\beta}(\Omega)\), for some \(\beta\in(0,1)\). We aim to show that the sequence\,\,\(\{u_n\}_n\) converges, up to a subsequence, to some \(u\in C^2(\Omega)\) solving\,\,$\mathcal{P}(\Omega,\phi)$. Observe that it is not at all clear that the \(C^{2,\beta}\) norms of the \(u_n\)'s will be uniformly bounded. Thus, to pass to the limit (and to ensure that the limit still solves \eqref{equation-weingarten}) we will need to prove H\"older estimates for the sequence \(\{u_n\}_n\). These rely on the fact that the Gauss map of the graph of a function solving \eqref{equation-weingarten}-\eqref{equation-uniformellipticity}-\eqref{equation-minimaltype} is a \textit{quasiregular mapping} (see \cite{simon-quasiconformal,galvezquasiconformal,galvezbarbieri2026asymptotic}).

Indeed, let \(u\in C^2(B_R)\) be a solution of \eqref{equation-weingarten}-\eqref{equation-uniformellipticity}-\eqref{equation-minimaltype}. Let \(K=k_1k_2\leq0\) be the gaussian curvature of the graph of \(u\). Then, uniform ellipticity implies:
\begin{equation}\label{equation-quasiregular}
	2|K|\leq k_1^2+k_2^2\leq\frac{2}{\Lambda}|K|,
\end{equation}
which implies quasiregularity of the Gauss map (see \cite{galvezbarbieri2026asymptotic}).

Now, observe that the definition of uniform ellipticity for Weingarten surfaces is, in general, different from the definition of uniform ellipticity for fully nonlinear elliptic equations (see \cite{caffarellicabre,galvezbarbieri2026asymptotic}, in particular Remark 2.4 in \cite{galvezbarbieri2026asymptotic}):
\begin{definition}\label{definition-uniformellipticity}
	Let \(\mathcal{S}^2\) be the set of \(2\times2\) symmetric matrices. An operator \(F:\R^2\times\mathcal{S}^2\to\R\) is called \textit{uniformly elliptic} if there exist positive constants \(\mu\), $\bar{\lambda}$, $\bar{\Lambda}$, with \(\bar{\lambda}\leq\bar{\Lambda}\), such that for any \(p,q\in\R^2\) and any \(M,N\in\mathcal{S}^2\), we have:
	\begin{equation}\label{equation-ellipticityF}
		F(p,M)-F(q,N)\leq\bar{\Lambda}|(M-N)^+|-\bar{\lambda}|(M-N)^-|+\mu |p-q|,
	\end{equation}
	where \((M-N)^+\) and \((M-N)^-\) are the positive and negative part respectively of the symmetric matrix \(M-N\). The constants \(\mu\), \(\bar{\lambda}\), \(\bar{\Lambda}\) are called \textit{ellipticity constants}.
\end{definition}
The following result will give the \(C^{2,\alpha}\) estimates we need.
\begin{lemmaletter}[\cite{galvezbarbieri2026asymptotic}, Lemma 2.5]\label{lemma-C2alphaestimate}
	Let \(R>0\) and suppose that \(u\in C^2(B_R)\) is a solution of \(F(Du,D^2u)=0\) in \(B_R\), with \(F\) uniformly elliptic. Then there exists \(\alpha\in(0,1)\), depending only on $\bar{\lambda}$, $\bar{\Lambda}$, such that \(u\in C^{2,\alpha}(\bar{B}_{\frac{R}{2}})\) and the following estimates hold:
	\begin{equation*}
		\sum_{k=1}^{2}R^k||D^ku||_{L^\infty (\bar{B}_{\frac{R}{2}})}+R^{2+\alpha}[D^2u]_{C^\alpha(\bar{B}_{\frac{R}{2}})}\leq C\left(||u||_{L^\infty (B_R)}+R^2|F(0,0)|\right),
	\end{equation*}
	where \(C\) depends only on $\bar{\lambda}$, $\bar{\Lambda}$ and \(R\mu\).
\end{lemmaletter}
Then, we only need to show that \(u\) solves a uniformly elliptic equation, for some\,\,\(F\). This comes as a consequence of \(C^2\) estimates for \(u\).
\begin{lemmaletter}[\cite{galvezbarbieri2026asymptotic}, Lemma 2.8]\label{lemma-Funiformelliptic}
	Let \(u\) define a solution of \eqref{equation-weingarten}-\eqref{equation-uniformellipticity} in a domain \(\Omega\subset\R^2\), where we only assume that \(f\in C^{0,1}_{\text{loc}}([c,+\infty))\) and \(f(c)=c\geq0\). Then, we can rewrite \eqref{equation-weingarten} as a solution to a fully nonlinear elliptic equation:
	\begin{equation}\label{equation-F}
		F(Du,D^2u)=0\quad\text{in }\Omega.
	\end{equation}
	Note that \(f\) may be not defined at 0. To ensure that \(F\) is well defined for all \((p,M)\in\R^2\times\mathcal{S}^2\), we extend the domain of \(f\) from \([c,+\infty)\) to \((-\infty,+\infty)\) by reflecting the graph of \(f\) with respect to the line \(k_2=k_1\) (i.e. \(f\circ f=\text{Id}\)).
	
	Moreover, if \(||u||_{C^2(\overline{\Omega})}\leq K\) for some \(K>0\) constant, then \(F\) can be enhanced to be uniformly elliptic, as in \eqref{equation-ellipticityF}.
\end{lemmaletter}
Now, the question is reduced to obtaining \(C^2\) estimates. These will come as a consequence of curvature estimates and \(C^1\) estimates.
\begin{lemmaletter}[\cite{galvezbarbieri2026asymptotic}, Lemma 2.6]\label{lemma-curvatureestimates}
	Let $\Sigma$ be a complete, uniformly elliptic Weingarten surface of minimal type, with boundary $\partial\Sigma$. Here, we only assume \(f\in C^{0,1}_{\text{loc}}[0,+\infty)\). Suppose that the image of its Gauss map lies in an open hemisphere of \(\esse^2\). Then:
	\begin{equation}\label{equation-curvatureestimate}
		|\sigma(p)|_g\leq\frac{C}{d(p,\partial\Sigma)}\quad\text{for all } p\in\Sigma,
	\end{equation}
	where \(C\) depends only on $\Lambda$. Here, \(d(\cdot,\cdot)\) is the intrinsic distance on $\Sigma$ and \(|\sigma(p)|_g\) denotes the intrinsic norm of its second fundamental form at \(p\in\Sigma\).
\end{lemmaletter}

From the quasiregularity of the Gauss map, gradient estimates follow:
\begin{theorem}\label{theorem-gradientestimate}
	Let \(u:B_R\to\R\) solve \eqref{equation-weingarten}-\eqref{equation-uniformellipticity}-\eqref{equation-minimaltype}. Then:
	\begin{equation*}
		|Du(0)|\leq C_1\exp\left\{C_2\frac{\text{osc}_{B_R}u}{R}\right\},
	\end{equation*}
	with \(\text{osc}_{B_R}u\defeq\sup_{B_R}u-\inf_{B_{R}}u\), where \(C_1\) and \(C_2\) only depend on $\Lambda$.
\end{theorem}
\begin{proof}
	The result follows directly from the results in \cite{simon-quasiconformal}. Indeed, by \eqref{equation-quasiregular}, the Gauss map of the graph of \(u\) is quasiregular, and we can use the results of \cite{simon-quasiconformal}, in particular \cite[Theorem 4.2]{simon-quasiconformal} and its Corollary.
\end{proof}
\begin{remark}
	Observe that, for any $\epsilon\in(0,1)$, Theorem \ref{theorem-gradientestimate} gives:
	\begin{equation*}
		||Du||_{L^\infty(B_{\epsilon R})}\leq C_1 \exp\left\{C_2\frac{||u||_{L^\infty(B_{R})}}{(1-\epsilon)R}\right\},
	\end{equation*}
	with \(C_1\) and \(C_2\) only depending on $\Lambda$.
\end{remark}
Once we have Lemma \ref{lemma-curvatureestimates} and Theorem \ref{theorem-gradientestimate}, \(C^2\) estimates follow easily. Indeed, consider \(u\in C^2(B_R)\) solution of \eqref{equation-weingarten}-\eqref{equation-uniformellipticity}-\eqref{equation-minimaltype}. Then:
\[D^2u=\sqrt{1+|D u|^2}\sigma,\]
which implies:
\[|D^2u|\leq(1+|Du|^2)^{\frac{3}{2}}|\sigma|_g,\]
where \(|\cdot|\) is the euclidean norm and \(|\cdot|_g\) is the norm with respect to the intrinsic metric of the graph of \(u\). Then, for any $\epsilon\in(0,1)$:
\begin{equation}\label{equation-C2estimate}
	||D^2u||_{L^\infty(B_{\epsilon R})}\leq\frac{C}{(1-\epsilon )R}\exp\left\{\tilde{C}\frac{||u||_{L^\infty(B_{R})}}{(1-\epsilon)R}\right\},
\end{equation}
where \(C\) and $\tilde{C}$ only depend on the ellipticity constant $\Lambda$. Last equation, jointly with Lemma \ref{lemma-Funiformelliptic} and Lemma \ref{lemma-C2alphaestimate}, gives the desired \(C^{2,\alpha}\) estimates in terms of \(C^0\) estimates for \(u\).

\begin{lemma}\label{lemma-boundedconvergence}
	Let $\Omega$ be an open bounded domain of \(\R^2\) and consider a sequence \(\{u_n\}_{n\in\N}\subset C^2(\Omega)\) solving \eqref{equation-weingarten}-\eqref{equation-uniformellipticity}-\eqref{equation-minimaltype}, such that:
	\[\sup_{n\in\N}||u_n||_{C^0(\Omega)}<+\infty.\]
	Then, \(\{u_n\}_n\) admits a subsequence converging in \(C^2_{\text{loc}}(\Omega)\) to a function solving \eqref{equation-weingarten}-\eqref{equation-uniformellipticity}-\eqref{equation-minimaltype}.
\end{lemma}
\begin{proof}
	Consider \(y\in\Omega\) and \(R>0\) such that \(\overline{B_R(y)}\subset\Omega\). Then, there exists $\delta>0$ such that \(B_{(1+\delta)R}(y)\subset\Omega\), and by \eqref{equation-C2estimate} we find:
	\[\sup_{n\in\N}||u_n||_{C^2(B_R(y))}<+\infty.\]
	Thus, Lemma \ref{lemma-Funiformelliptic} implies that the \(u_n\)'s solve a uniformly elliptic equation, and Lemma \ref{lemma-C2alphaestimate} provides uniform \(C^{2,\alpha}\) estimates for the sequence. By Ascoli-Arzelà, we then get that a subsequence of \(\{u_n\}_n\) converges to some \(C^2\) function \(u\) solving \eqref{equation-weingarten}-\eqref{equation-uniformellipticity}-\eqref{equation-minimaltype}.
	
	Covering $\Omega$ with a countable number of open discs and by a diagonal trick, we see that actually \(u\) can be extended to the whole $\Omega$, \(u\) solves \eqref{equation-weingarten}-\eqref{equation-uniformellipticity}-\eqref{equation-minimaltype} and there is a subsequence of \(\{u_n\}_n\) converging to \(u\) in \(C^2_\text{loc}(\Omega)\).
\end{proof}
We are now ready to solve the Dirichlet problem $\mathcal{P}(\Omega,\phi)$ for $\phi\in C^0(\partial\Omega)$.
\begin{theorem}\label{theorem-dirichletC0}
	Let $\Omega \subset \mathbb{R}^2$ be a strictly convex, bounded domain of class
	$C^{2,\alpha}$, with
	$\alpha\in(0,1)$, and let $\varphi \in C^{0}(\partial\Omega)$. Then there exists a unique solution
	$u \in C^{2}(\Omega)\cap C^0(\overline{\Omega})$ to \eqref{equation-weingarten}-\eqref{equation-uniformellipticity}-\eqref{equation-minimaltype}, where
	$f \in C_{\mathrm{loc}}^{0,1}[0,+\infty)$, subject to the boundary condition:
	\[
	u = \varphi \qquad \text{on } \partial\Omega.
	\]
\end{theorem}
\begin{proof}
	As anticipated, let $\{\phi_n\}_n\subset C^\infty(\partial\Omega)$ be a sequence such that:
	\[\phi_n\xrightarrow{C^0}\phi.\]
	For \(n\in\N\), let \(u_n\) be the unique solution to $\mathcal{P}(\Omega,\phi_n)$. It is provided by Theorem\,\,\ref{theorem-galvezdirichlet}. Observe that, by the maximum principle, we have that \(\{u_n\}_n\) is bounded in \(C^0(\Omega)\), so that Lemma \ref{lemma-boundedconvergence} can be applied: there exists \(u\in C^2(\Omega)\) solving \eqref{equation-weingarten}-\eqref{equation-uniformellipticity}-\eqref{equation-minimaltype} and a subsequence of \(\{u_n\}_n\), that we do not relabel, converging to \(u\) in \(C^2_\text{loc}(\Omega)\).
	
	We want to show that \(u\) solves \(\mathcal{P}(\Omega,\phi)\). But this is immediate since, by the maximum principle for \eqref{equation-weingarten}-\eqref{equation-uniformellipticity}-\eqref{equation-minimaltype}, we get:
	\begin{equation*}
	||u_k-u_j||_{C^0(\overline{\Omega})}\leq||\phi_k-\phi_j||_{C^0(\partial\Omega)},
	\end{equation*}
	which implies that \(\{u_n\}_n\) is a Cauchy sequence in \(C^0(\overline{\Omega})\). Then, \(u\in C^0(\overline{\Omega})\) and:
	\[u_n\xrightarrow{C^0(\overline{\Omega})}u.\]
	In particular: \[u=\phi\text{ on }\partial\Omega.\] Uniqueness comes immediately from the maximum principle.
\end{proof}
	\section{The Perron method and geometric barriers}\label{section-perron}
	Let $\Omega$ be a bounded open domain of \(\R^2\), and consider a continuous function \(u:\overline{\Omega}\to\R\). Given a small closed disc \(B\subset\Omega\), let $\tilde{u}_B$ be the unique function solving \eqref{equation-weingarten}-\eqref{equation-uniformellipticity}-\eqref{equation-minimaltype} on \(B\) with boundary value:
	\[\tilde{u}_B|_{\partial B}=u|_{\partial B}.\]
	That is, $\tilde{u}_B$ is the solution to $\mathcal{P}(B,u|_{\partial B})$. We know that $\tilde{u}_B$ is well defined by Theorem \ref{theorem-dirichletC0}.
	
 We say that \(u\) is a \textit{subsolution} to \(\mathcal{P}(\Omega,\phi)\) (resp. \textit{supersolution}) if, for any closed disc \(B\subset\Omega\), we have:
	\begin{equation*}
		u|_B\leq\tilde{u}_B\quad\text{(resp. }u|_B\geq\tilde{u}_B\text{)},
	\end{equation*}
	and:
	\begin{equation*}
		u|_{\partial\Omega}\leq \phi\quad\text{(resp. }u|_{\partial\Omega}\geq \phi\text{)}.
	\end{equation*}
	We also define the concept of \textit{lifting}:
	\begin{definition}
		Let $\Omega$ be a bounded domain in \(\R^2\), \(u\in C^0(\Omega)\), \(B\subset\Omega\) a small closed disc. We define the \textit{lifting} of \(u\) with respect to the disc \(B\) as:
		\begin{equation*}
			U_B(x)=
			\begin{cases}
				\tilde{u}_B(x), & x\in B,\\
				u(x),       & x\in \Omega\setminus B.
			\end{cases}.
		\end{equation*}
	\end{definition}
The following properties follow easily:
\begin{enumerate}
	\item For a function \(u\in C^2(\Omega)\cap C^0(\overline{\Omega})\) such that \(u|_{\partial\Omega}\leq\phi\) (resp. \(u|_{\partial\Omega}\geq\phi\)), being a subsolution (resp. a supersolution) is equivalent to \(k_2\geq f(k_1)\) (resp. \(k_2\leq f(k_1)\)) on the graph of \(u\);
	\item The supremum of a finite family of subsolutions is again a subsolution;
	\item The infimum of a finite family of supersolutions is again a supersolution;
	\item For any \(\overline{u}\) supersolution of $\mathcal{P}(\Omega,\phi)$ and any \(\underline{u}\) subsolution of $\mathcal{P}(\Omega,\phi)$, we have $\underline{u}\leq\overline{u}$;
	\item Any lifting of a subsolution is still a subsolution.
\end{enumerate}
\begin{definition}
	Let $\Omega$ be a bounded domain. We say that \(p\in\partial\Omega\) admits a \textit{barrier} for \(\mathcal{P}(\Omega,\phi)\) if there exist a subsolution \(\underline{\psi}\) and supersolution \(\overline{\psi}\) both in \(C^0(\overline{\Omega})\) such that \(\underline{\psi}(p)=\phi(p)=\overline{\psi}(p)\).
\end{definition}

\begin{proposition}[Perron Process]\label{proposition-perron}
	Given a bounded domain \(\Omega\subset\R^2\), let $\phi\in C^0(\partial\Omega)$, and define:
	\begin{equation*}
		\mathcal{S}\defeq\left\{v\,|\,v\text{ is a subsolution of }\mathcal{P}(\Omega,\phi)\right\}.
	\end{equation*}
	Assume that $\mathcal{S}\neq\emptyset$. Then:
	\begin{itemize}
		\item[(i)] The function \(u(x)\defeq \sup_{v\in\mathcal{S}}v(x)\) is in \(C^2(\Omega)\) and satisfies \eqref{equation-weingarten};
		\item[(ii)] Suppose that \(p\in\partial\Omega\) admits a barrier for $\mathcal{P}(\Omega,\phi)$. Then \(u\) extends continuously to \(p\), by setting \(u(p)=\phi(p)\).
	\end{itemize}
\end{proposition}
\begin{proof}
	Once Lemma \ref{lemma-boundedconvergence} and Theorem \ref{theorem-dirichletC0} are settled, the proof of (i) is analogous to the proof of the Perron process for the Laplace equation, which can be found in \cite[Theorem 2.12]{gilbargtrudinger}.
	
	Regarding (ii), the result is straightforward since, by the properties of sub- and supersolutions, we have:
	\[\underline{\psi}\leq u\leq\overline{\psi},\]
	and \(\underline{\psi}\), $\overline{\psi}$ are continuous up to the boundary, with:
	\[\underline{\psi}(p)=\phi(p)=\overline{\psi}(p).\]
\end{proof}
The barriers for the Perron process will be built using the rotational examples of minimal type Weingarten surfaces. These have been classified in \cite{fernandezweingarten}.

Fix a function \(f\in C^{1,1}_{\text{loc}}[0,+\infty)\), with \(f(0)=0\), satisfying \eqref{equation-uniformellipticity}, and consider the Weingarten equation \eqref{equation-weingarten}. According to \cite{fernandezweingarten, galvezbarbieri2026asymptotic}, for any \(R>0\), there exists a rotational solution to \eqref{equation-weingarten} on \(\overline{B_R(0)}^C\) with zero boundary value, whose behavior at infinity depends on \(f'(0)\):

\begin{lemmaletter}[\cite{galvezbarbieri2026asymptotic}, Lemma 6.1]\label{lemma-galvez-continuousfamily}
	Let $f\in C^{1,1}_{loc}[0,+\infty)$ satisfy \eqref{equation-uniformellipticity}-\eqref{equation-minimaltype}. For any $R_0>0$ and $C_0> 0$, there exists a unique radial solution $u$ of
	\begin{equation}\label{e8.0}
		\kappa_2=f(\kappa_1)\quad\mbox{in}~~\overline{B_{R_0}}^C
	\end{equation}
	with the initial conditions $u(R_0)=0$ and $u'(R_0)=C_0$. In addition, $u$ has the following properties:~\\[1mm]
	(i) The solution $u$ is nondecreasing, $u'$ is nonincreasing and $u'(r)\to 0$ as $r\to \infty$. Moreover,
	\begin{equation*}\label{e8.0-2}
		0\leq u(r)\leq \frac{C_0R_0^{\Lambda}}{1-\Lambda}r^{1-\Lambda}, \quad  0<u'(r)\leq C_0, \quad 
		-\frac{C_0(1+C_0^2)}{\Lambda R_0}\leq u''(r)<0,
		~\forall ~r\in [R_0,+\infty).
	\end{equation*}
	(ii) If $-1<f'(0)<0$, there exists $c> 0$ such that
	\begin{equation}\label{e8.1}
		\lim_{r\to \infty} \frac{u(r)}{r^{1+f'(0)}}=c.
	\end{equation}
	Moreover, $c$ is a continuous function of $R_0,C_0$ and
	\begin{equation}\label{e8.2-1}
		c\to 0 ~~\mbox{ as }~~R_0\to 0 ~~(\mbox{or }~~C_0\to 0)\quad\mbox{and}\quad c\to \infty ~~\mbox{ as }~~R_0\to \infty.
	\end{equation}
	(iii) If \(f'(0)=-1\), there exists \(c>0\) such that
		\begin{equation}\label{e8.1a}
		\lim_{r\to \infty} \frac{u(r)}{\log r}=c.
	\end{equation}
	Moreover, $c$ is a continuous function of $R_0,C_0$ and
	\begin{equation}\label{e8.2-1a}
		c\to 0 ~~\mbox{ as }~~R_0\to 0 ~~(\mbox{or }~~C_0\to 0)\quad\mbox{and}\quad c\to \infty ~~\mbox{ as }~~R_0\to \infty.
	\end{equation}
	(iv) If $f'(0)<-1$, there exists $h> 0$ such that
	\begin{equation*}\label{e8.1-2}
		\lim_{r\to \infty} u(r)=h.
	\end{equation*}
	Moreover, $h$ is a continuous function of $R_0,C_0$ and
	\begin{equation}\label{e8.2-2}
		h\to 0 ~~\mbox{ as }~~R_0\to 0 ~~(\mbox{or }~~C_0\to 0)\quad\mbox{and}\quad h\to \infty ~~\mbox{ as }~~R_0\to \infty.
	\end{equation}
\end{lemmaletter}
\begin{remark}
	Observe that points (i), (ii), (iv) are proved in \cite[Lemma 6.1]{galvezbarbieri2026asymptotic}. The logarithmic case of (iii) is not stated explicitly, but it follows easily from the same argument.
\end{remark}
\section{On the exterior Dirichlet problem}\label{section-exteriordirichlet}
In this section, we solve the Dirichlet problem on the exterior of a domain $\Omega$, with zero boundary value, assuming some condition on $\Omega$. In particular, we show the existence of non rotational Weingarten annular ends. Let $\Omega$ be a bounded domain and assume that there exists a \(R>0\) such that, for each point \(p\in\partial\Omega\) there exists a point \(q_p\in\Omega\) for which \(B_R(q_p)\subset\Omega\) and $p\in\partial B_R(q_p)$. In this case, we say that $\Omega$ satisfies an \textit{interior radius }\(R\) \textit{sphere condition}.
\begin{theorem}
	Let $\Omega\subset\R^2$ be a bounded domain, such that $\overline{\Omega}^C$ is connected, satisfying an interior radius \(R\) sphere condition. Then, there exists a non trivial solution \(u\in C^2(\overline{\Omega}^C)\cap C^0(\Omega^C)\) to $\mathcal{P}(\overline{\Omega}^C,0)$, with limit unit normal \(e_3\). 
\end{theorem}
\begin{proof}
	Let \(n\in\N\) be sufficiently big so that \(\overline{\Omega}\subset B_n(0)\). We want to apply the Perron process to the domains:
	\[A_n\defeq B_n(0)\setminus\overline{\Omega}.\]
	They form an exhaustion of $\overline{\Omega}^C$ and so we aim to construct a solution to $\mathcal{P}(\overline{\Omega}^C,0)$ by taking the limit of the solutions found on the bounded domains $A_n$'s. The problem we want to solve on \(A_n\) is $\mathcal{P}(A_n,\phi_n)$, with:
	\[\phi_n(x)=
	\begin{cases}
		0\quad&\text{if }x\in\partial\Omega;
		\\
		L_n&\text{if }x\in\partial B_n(0).
	\end{cases}\]
	Here \(L_n\) is a constant chosen as follows. By Lemma \ref{lemma-galvez-continuousfamily}, we know that there is a solution \(v\) to $\mathcal{P}(\overline{B_R(0)}^C,0)$. Then \(L_n\) is chosen so that, for each \(p\in\partial\Omega\), one has:
	\[v(x-q_p)\geq L_n\quad\text{on }\partial B_n(0).\]
	Remember that \(q_p\in\Omega\) is the center of a disc satisfying the interior radius \(R\) sphere condition on \(p\). By compactness of $\partial\Omega\cup\partial B_n(0)$, and the asymptotic of \(v\), a choice of the \(L_n\)'s is possible so that:
	\[L_n= v|_{\partial B_n(0)}+o(1)\text{ as }n\to+\infty.\]
	To show that \(\mathcal{P}(A_n,\phi_n)\) is solvable, we need to show that it admits a barrier at each point of $\partial A_n$. At a point $p\in\partial \Omega$, the inferior barrier will be the identically zero function, while the superior barrier will be:
	\[v(x-q_p)|_{A_n},\]
	which is a supersolution by the definition of \(L_n\). At every point \(p\in\partial B_n(0)\), convexity permits us to use a horizontal plane and a tilted plane as upper and lower barrier respectively. Then, $\mathcal{P}(A_n,\phi_n)$ admits a unique solution \(v_n\in C^2(A_n)\cap C^0(\overline{A_n})\). Now, the result of Section \ref{section-dirichlet}, together with a diagonal argument, can be used to extract a subsequence from \(\{v_n\}_{n\in\N}\), converging in \(C^2_\text{loc}(\overline{\Omega}^C)\) to a function \(u\) solving the Weingarten equation. The existence of barriers at points of $\partial\Omega$ also imply that \(u\in C^2(\overline{\Omega}^C)\cap C^0(\Omega^C)\), with \(u|_{\partial\Omega}\equiv0\). It remains to show that \(u\) is non trivial. To prevent this, observe the following. Let \(p\in\partial\Omega\), and consider the point \(q_p\in\Omega\). For \(K\defeq\max_{x\in\partial\Omega}v(x-q_p)\), consider the function \(\left(v(x-q_p)-K-\delta\right)|_{\overline{\Omega}^C}\), for some $\delta>0$, which is negative on $\partial\Omega$. By the choice of the \(L_n\)'s and the maximum principle, we have that, for \(n\) sufficiently big:
	\[(v-K-\delta)|_{A_n}\leq v_n.\]
	Then, this property is preserved in the limit: \(v(x)-K-\delta\leq u(x)\) for \(|x|\) sufficiently big. In the case in which \(f'(0)\geq-1\), this automatically implies that \(u\) is nontrivial, since \(v(x)\to+\infty\) as \(|x|\to+\infty\). Moreover, as \(u\) is also below a vertical translate of \(v\), we have that \(u(x)=o(|x|)\), and then Theorem \ref{theorem-gradientestimate} implies that \(|Du(x)|\to0\) as \(|x|\to+\infty\), i.e. the graph of \(u\) has limit unit normal \(e_3\). When \(f'(0)<-1\), the standard maximum principle and the choice of \(K\) imply that, for $\delta>0$ sufficiently small, there exists an $\epsilon>0$ such that \(v(x)-K-\delta\to\epsilon\) as \(|x|\to+\infty\), thus:
	\[\liminf_{|x|\to+\infty}u(x)\geq\epsilon>0.\]
	Also in this case, proving that the graph of the solution has limit unit normal \(e_3\) is straightforward.
\end{proof}
\section{Maximum principle at infinity}\label{section-maximumprincipleinfinity}
Barbieri, Gálvez, Lian and Zhang described the asymptotic behavior of exterior solutions to the uniformly elliptic Weingarten equation of minimal type. In particular, they showed that the behavior is analogous to the one of rotational solutions, which had been classified by Fernández and Mira in \cite{fernandezweingarten}:
\begin{theoremletter}[\cite{galvezbarbieri2026asymptotic}, Theorem 1.3]\label{theorem-galvez-c}
	Let $u$ define a uniformly elliptic Weingarten graph of minimal type over $\overline{B_R}^C$ with $u$ bounded from below and limit unit normal \(e_3\). Then
	\begin{equation}\label{e1.1}
		\left\{
		\begin{aligned}
			&u(x)\simeq |x|^{1+f'(0)}~~
			&\mbox{or}\quad u_{\infty}\in \mathbb{R},~~\quad&\mbox{if}~~-1<f'(0)<0;\\
			&u(x)\simeq \log |x|~~
			&\mbox{or}\quad u_{\infty}\in \mathbb{R},~~\quad&\mbox{if}~~f'(0)=-1;\\
			&u_{\infty}\in \mathbb{R},\quad~&~&\mbox{if}~~f'(0)<-1,\\
		\end{aligned}
		\right.
	\end{equation}
	where $u\simeq v$ means that there exists a positive constant $c$ such that
	\begin{equation*}
		\lim_{x\to \infty}\frac{u(x)}{v(x)}=c.
	\end{equation*}
	
	Moreover, if $f'(0)=-1$, there exist $d\geq 0$ and $c\in \mathbb{R}$ such that for any $0<\alpha<1$,
	\begin{equation}\label{e1.10}
		u(x)=d\log |x|+c+O(|x|^{-\alpha}),~\forall ~x\in \overline{B_R}^C
	\end{equation}
	and
	\begin{equation*}\label{e1.11}
		|Du(x)|= O(|x|^{-1}), \quad
		|D^2u(x)|= O(|x|^{-2}),~\forall ~x\in \overline{B_R}^C.
	\end{equation*}
	Here, $f=O(g)$ means that there exists a constant $C$ such that $|f|\leq Cg$.
\end{theoremletter}
Observe that, in the case in which \(f'(0)=-1\), they find the asymptotic expansion \eqref{e1.10}, which gives more information than the general asymptotic description \eqref{e1.1}. In fact, when \(f'(0)=-1\), we have differentiability of the associated equation: \[F(Du,D^2u)=0\]
at umbilic points, and therefore general results on fully nonlinear uniformly elliptic equations can be used (see \cite[Theorem 5.1]{galvezbarbieri2026asymptotic}).

Using the asymptotic expansion \eqref{e1.10}, they are able to prove a maximum principle at infinity in the case \(f'(0)=-1\) \cite[Corollary 1.5]{galvezbarbieri2026asymptotic}.
In this section, we carry out an asymptotic expansion, analogous to \eqref{e1.10}, in the case \(f'(0)\in (-1,0)\), for ends of unbounded height. This permits us to completely identify which classes of Weingarten minimal ends admit a maximum principle at infinity and study how the symmetries of the boundary datum are inherited by the solutions to the exterior Dirichlet problem. We first make the following remark: 

\begin{remark}\label{remark-unboundedends}
	Observe that the two classes \(f'(0)\in(-1,0)\) and \(f'(0)<-1\) represent the same class of unoriented surfaces. The different kind of behavior that seems to appear in Theorem \ref{theorem-galvez-c} depends only on the fact that, when we look at graphs, we fix an orientation. Moreover, we are only looking at ends which are bounded from below. For instance, the class \(f'(0)<-1\) admit ends which are unbounded from below, which have the same asymptotic as the ones described by Theorem \ref{theorem-galvez-c} in the class \(f'(0)\in(-1,0)\) (but with a minus sign for the coefficient \(c\)). This is because, if \(u\) solves the Weingarten equation \(k_2=f(k_1)\), with the orientation given by the upper normal, then \(-u\) solves \(k_2=\tilde{f}(k_1)\), with $\tilde{f}(t)\defeq-f^{-1}(-t)$. In particular, $\tilde{f}'(0)=f'(0)^{-1}$.
\end{remark}
The key to extend the asymptotic expansion is the fact that by Theorem \ref{theorem-galvez-c} we identify a coefficient \(c\), which guides the asymptotic of the Weingarten end of minimal type:
 \begin{definition}\label{definition-leadingcoefficient}
	Let \(u\) be a solution to $\mathcal{P}(\overline{\Omega}^C,\phi)$, bounded from below and with limit unit normal \(e_3\). Define the \textit{limiting height} of \(u\) as:
	\[\lim_{|x|\to+\infty}u(x)\in\R\cup\{+\infty\},\]and the \textit{leading coefficient} \(c\) of \(u\) as:
	\begin{itemize}
		\item \(c\defeq\lim_{|x|\to+\infty}\frac{u(x)}{\log|x|}\) if \(f'(0)=-1\);
		\item \(c\defeq\lim_{|x|\to+\infty}\frac{u(x)}{|x|^{1+f'(0)}}\) if \(-1<f'(0)<0\).
	\end{itemize} Observe that this definition is well posed by Theorem \ref{theorem-galvez-c} and Remark \ref{remark-unboundedends}. Moreover, \(c\) is finite, and \(c=0\) if and only if \(u\) is bounded.
\end{definition}
In principle, it is not at all clear that different solutions to the same Dirichlet problem must have different leading coefficients (or, in the case of bounded ends, different limiting heights). This will be a consequence of the comparison principle below (Theorem \ref{theorem-weakmaximumprincipleatinfinity}).

 Let $\Omega$ be a bounded smooth domain and consider \(u\in C^2(\overline{\Omega}^C)\cap C^0(\Omega^C)\) solution to \eqref{equation-weingarten}-\eqref{equation-uniformellipticity}-\eqref{equation-minimaltype} in \(\overline{\Omega}^C\), with \(f'(0)\in(-1,0)\). Assume that \(u\) is unbounded and let \(c>0\) be its leading coefficient. In the following, it will be convenient to use the polar orthonormal frame \(\{e_r,e_\theta\}\). Define the matrix \(Q\defeq(e_r,e_\theta)\).

\begin{lemma}\label{lemma-asymptotics-uv}
	For \(u\) as above, we have that:
	\begin{equation*}
		Du=c(1+f'(0))r^{f'(0)}e_r+o(r^{f'(0)})
	\end{equation*}
	and:
	\begin{equation*}
		Q^TD^2uQ=c(1+f'(0))r^{f'(0)-1}\begin{pmatrix}
			f'(0)&0\\
			0&1
		\end{pmatrix} +o(r^{f'(0)-1}),
	\end{equation*}
	where \(c\) is the leading coefficient of \(u\).
\end{lemma}
\begin{proof}
	Consider the function:
	\[u_R(y)\defeq\frac{u(Ry)}{R},\]
	which, for \(R>0\) sufficiently big, is defined on the annulus:
	\[A'\defeq B_4\setminus\overline{B_{\frac{1}{2}}}.\]
	Observe that \(u_R\) solves:
	\begin{equation}\label{equation-blowdown}
		k_2=f_R(k_1),
	\end{equation}
	for \(f_R(t)\defeq Rf(\frac{t}{R})\), and \eqref{equation-blowdown} is uniformly elliptic, with ellipticity constant\,\,$\Lambda$ independent of \(R\).
	
	We want to use the machinery of Section \ref{section-dirichlet} to prove that the family \(\{u_R\}_R\) is \(C^{2,\alpha}\)-bounded in some smaller annulus.
	
	Observe, nonetheless, that each \(u_R\) solves its own fully nonlinear equation as in:
	\begin{equation}\label{equation-FR}
		F_R(Du_R,D^2u_R)=0.
	\end{equation}
	Now, the \(C^{2,\alpha}\) bound provided by the lemmata of Section \ref{section-dirichlet} depends on the ellipticity constants of the fully nonlinear equation, which are given, as in Lemma \ref{lemma-Funiformelliptic}, by enhancing \eqref{equation-FR} to a uniformly elliptic equation. By reading through the proof of Lemma \ref{lemma-Funiformelliptic}, which can be found in \cite[Lemma 2.8]{galvezbarbieri2026asymptotic}, we see that every \(u_R\) solves a fully nonlinear, uniformly elliptic equation, whose constants of ellipticity depend only on:
	\begin{itemize}
		\item The ellipticity constant of \(f_R\) (which is uniform in \(R\)!);
		\item The \(C^2\) norm of \(u_R\).
	\end{itemize}
	We then need to bound the latter uniformly. But this is immediate, since the results coming from \cite{simon-quasiconformal} only use the quasiregularity of the Gauss map, which is preserved by the blow-down: in fact the bounds of Lemma \ref{lemma-curvatureestimates} and Theorem \ref{theorem-gradientestimate} only depend on the ellipticity constant of \(f_R\), which is equal to $\Lambda$.
	
	Then, we see that the second fundamental form of the graphs of the \(u_R\) is uniformly bounded in \(A'\). Moreover, the asymptotics of Theorem \ref{theorem-galvez-c} imply that \(u_R\to0\) uniformly in \(A'\), as \(R\to+\infty\), so that, in the smaller annulus
	\[A''\defeq B_3\setminus\overline{B_{\frac3{4}}}\subset A',\]
	we have:
	\[\sup_{R}||u_R||_{C^2(A'')}<+\infty.\]
	Then, by the proof of Lemma \ref{lemma-Funiformelliptic}, each \(u_R\) solves a fully nonlinear, uniformly elliptic equation, whose ellipticity constants are uniform in \(R\). By Lemma \ref{lemma-C2alphaestimate}, we get that there exists an $\alpha\in(0,1)$ and a constant \(C\) that does not depend on \(R\) such that, on the smaller annulus:
	\[A\defeq B_2\setminus\overline{B_1},\]
	the following estimate holds:
	\begin{equation}\label{equation-estimate-uR}
		||u_R||_{C^{2,\alpha}(A)}\leq C||u_R||_{C^0(A'')}\leq CR^{f'(0)},
	\end{equation}
	where the last inequality comes from Theorem \ref{theorem-galvez-c}. Define now:
	\[U_R(y)\defeq R^{-f'(0)}u_R(y)=\frac{u(Ry)}{R^{1+f'(0)}}.\]
	We already know, by Theorem \ref{theorem-galvez-c}, that:
	\[U_R(y)\to c|y|^{1+f'(0)}\quad\text{uniformly in }A.\]
	By \eqref{equation-estimate-uR} we have, moreover, that:
	\begin{equation*}
		\sup_{R}||U_R||_{C^{2,\alpha}(A)}<+\infty,
	\end{equation*}
	so that, by Ascoli-Arzelà, we can upgrade the convergence of \(\{U_R\}_R\) to:
	\begin{equation*}
		U_R\xrightarrow{C^2(A)}c|y|^{1+f'(0)},
	\end{equation*}
	which proves the thesis.
\end{proof}
Next, we prove that we can actually linearize the Weingarten equation sufficiently far out. Remember that:
\begin{equation}\label{equation-weingarten2}
	F(Du,D^2u)=H-\sqrt{H^2-K}-f\left(H+\sqrt{H^2-K}\right),
\end{equation}
where $H$ and $K$ are the mean and gaussian curvature respectively of the graph of \(u\). We remark that, for \(f'(0)=-1\), the operator \(F\) is locally of class \(C^{1,1}\) in the region where \(K\leq0\) -- which is the case for the solutions we are considering -- which permits to linearize \eqref{equation-weingarten2}. In our case, we don't have such regularity for \(F\), therefore we need to take care of the umbilical locus, which is exactly where regularity might fail. Let \(v\) be another solution with the same leading coefficient as \(u\). Then:
\begin{lemma}\label{lemma-differentiability}
	With \(u\), \(v\) as above, define:
	\begin{equation*}
		p_t(y)\defeq tDu(y)+(1-t)Dv(y),\quad\quad M_t(y)\defeq tD^2u(y)+(1-t)D^2v(y).
	\end{equation*}
	Then, there exists \(R_0>0\) such that, for \(y\in B_{R_0}^C\), the map \(t\in[0,1]\to F(p_t(y),M_t(y))\) is differentiable. In particular, the graphs of \(u\) and \(v\) do not have umbilical points in \(B_{R_0}^C\).
\end{lemma}
\begin{proof}
	The operator \(F\) is of class \(C^1\) in the region in which \(K<0\). Since our solutions satisfy \(K\leq0\), we only need to take care of the umbilical locus. That is, we only need to prove that, for \(|y|\) sufficiently big, there are no umbilical points on the segment \((p_t,M_t)\). Remember that:
	\begin{equation}\label{equation-HK}
		H(p,M)=\frac{W^2\tr M-p^TMp}{2W^3},\quad K(p,M)=\frac{\det M}{W^4}, 
	\end{equation}
	with \(W(p)\defeq\sqrt{1+|p|^2}\), and observe that the estimates  of Lemma \ref{lemma-asymptotics-uv} remains valid on the whole segment \((p_t,M_t)\). Define for simplicity:
	\[\rho\defeq c(1+f'(0))r^{f'(0)-1}.\]Then:
	\begin{equation}\label{equation-estimateHK}
		H(p_t,M_t)=\rho\left(\frac{1+f'(0)}{2}+o(1)\right),\quad K(p_t,M_t)=\rho^2\left(f'(0)+o(1)\right),
	\end{equation}
	therefore:
	\[H^2(p_t,M_t)-K(p_t,M_t)=\rho^2\left(\frac{(1-f'(0))^2}{4}+o(1)\right),\]
	which is positive for \(|y|\) sufficiently big.
\end{proof}

Now, pick the \(R_0>0\) given by Lemma \ref{lemma-differentiability}. Then, the Weingarten equation \eqref{equation-weingarten2} can be differentiated, in \(B_{R_0}^C\), along the segment joining \(u\) to \(v\). Define:
\[w\defeq u-v.\]
The function \(w\) solves:
\begin{equation}\label{equation-differentiated}
	Lw=a^{ij}w_{ij}+b^iw_i=0\quad\text{in }\overline{B_{R_0}}^C,
\end{equation}
with \(a^{ij}(x)\) and \(b^i(x)\) given by:
\[a^{ij}(x)=\int_0^1(F_{M})_{ij}(p_t(x),M_t(x))dt,\quad b^i(x)=\int_0^1F_{p_i}(p_t(x),M_t(x))dt.\]
We use the following convention. Given a differentiable function $\Phi$ defined on the space of symmetric matrices, $\Phi_M$ will be the symmetric matrix representing the differential of $\Phi$, that is:
\[d\Phi=(\Phi_M)_{ij}dM_{ij}.\]
Since we can assume \eqref{equation-weingarten2} to be uniformly elliptic, we have that also \eqref{equation-differentiated} is uniformly elliptic. In fact, it is easy to show that the ellipticity constants of \(L\) are \(\frac{(1+o(1))\Lambda}{(1+|p|^2)^{3/2}}\) and \(\frac{1+o(1)}{\Lambda\sqrt{1+|p|^2}}\). We now study the asymptotic behavior of the coefficients of \(L\):
\begin{lemma}
	The coefficients \(A\defeq \{a^{ij}\}_{i,j}\) and \(b\defeq\{b^i\}_i\) of the linearized operator \(L\) satisfy:
	\begin{equation}\label{equation-estimateA}
		A=e_re_r^T-f'(0)e_\theta e_\theta^T + o(1),
	\end{equation}
	and:
	\begin{equation}\label{equation-estimateb}
		|b|=O(r^{2f'(0)-1}).
	\end{equation}
	
\end{lemma}
\begin{proof}
	We need to compute the derivatives of \(F\) along the segment \((p_t,M_t)\), with \(p_t\), \(M_t\) defined as in Lemma \ref{lemma-differentiability}. By direct calculations we see that, on points of differentiability of \(F\), one has:
	\begin{equation}\label{equation-101126}
		F_w=(1-f')H_w-\frac1{2}(1+f')\frac{(H^2-K)_w}{\sqrt{H^2-K}},
	\end{equation}
	for each \(w\in\{p_1,p_2,M_{11},M_{12},M_{22}\}\). Observe that \(f'\) is evaluated at \(k_1=H+\sqrt{H^2-K}\), and therefore \(f'=f'(0)+o(1)\). Let us start by computing the derivatives of \(H\) and \(K\) with respect to the entries of \(M\).
	Remember the formulas \eqref{equation-HK} for \(H\) and \(K\). One has:
\begin{equation}\label{equation-HKderivative}
		H_M=\frac{W^2I-pp^T}{2W^3},\quad K_M=\frac{\text{cof}M}{W^4}.
\end{equation}
Combining this formula with \eqref{equation-estimateHK} and the estimates of Lemma \ref{lemma-asymptotics-uv}, we get:
	\begin{equation*}
		Q^TH_MQ=\frac1{2}I+o(1),\quad Q^TK_MQ=\rho\begin{pmatrix}
			1&0\\
			0&f'(0)
		\end{pmatrix} +o(r^{f'(0)-1})
	\end{equation*}
	Then:
	\begin{equation*}
		Q^T\frac{(H^2-K)_M}{\sqrt{H^2-K}}Q=\begin{pmatrix}
			-1&0\\
			0&1
		\end{pmatrix}+o(1),
	\end{equation*}
	which implies:
	\begin{equation*}
		Q^TF_MQ=\begin{pmatrix}
			1&0\\
			0&-f'(0)
		\end{pmatrix}+o(1).
	\end{equation*}
	Now, rotating back to the euclidean coordinates, we find:
	\begin{equation*}
		F_M(p_t,M_t)=e_re_r^T-f'(0)e_\theta e_\theta^T+o(1),
	\end{equation*}
	with the \(o(1)\) going to zero, uniformly in \(t\), as \(|x|\) grows. Integrating on the segment \([0,1]\) we then find \eqref{equation-estimateA}.
	
	We now bound the first order terms. By direct computation, we see that:
	\[
	F_p(p,M)
	=-\frac{\operatorname{tr}(F_M M)\,p+2F_MMp}{1+|p|^2}.
	\]
	This expression, with the bound on the ellipticity constants of \(L\), gives:
	\[|F_p(p,M)|\leq\frac{2|F_MM||p|+2|F_MMp|}{1+|p|^2}\leq\frac{4|M||p|}{\Lambda (1+|p|^2)^{\frac3{2}}}.\]
	Now, remember that:
	\[|Du|^2=u_r^2+\frac{u_\theta^2}{r^2}\]
	and:
	\[|D^2u|^2=u_{rr}^2+2\left(\frac{u_{r\theta}}{r}-\frac{u_\theta}{r^2}\right)^2+\left(\frac{u_{\theta\theta}}{r^2}+\frac{u_r}{r}\right)^2,\]
	so that, by the estimates of Lemma \ref{lemma-asymptotics-uv}, we get:
	\[|Du|=O(r^{f'(0)})\quad\text{and}\quad |D^2u|=O(r^{f'(0)-1}).\]
	Observe that the same estimates hold on the whole segment joining the two solutions \(u\) and \(v\). Then:
	\[|b_i|=O(r^{2f'(0)-1}).\]
\end{proof}
We are now ready to give the following asymptotic expansion:
\begin{theorem}\label{theorem-asymptotic}
	Let \(u:B_R^C\to\R\) define a uniformly elliptic Weingarten graph of minimal type, with \(f'(0)\in(-1,0)\), and assume that \(u\) is unbounded, with limit unit normal \(e_3\). Then, given \(v\) a rotational solution with the same leading coefficient as \(u\), there exists two constants $\alpha\in(0,1)$, \(b\in\R\) such that:
	\begin{equation}\label{equation-asymptoticNEW}
		u(x)=v(|x|)+b+O(|x|^{-\alpha}).
	\end{equation}
	
\end{theorem}
\begin{proof}
	We first show that \(u-v\) is bounded. First of all observe that, by choosing \(\tilde R\) larger if necessary -- which does not change the study of the asymptotic expansion -- Lemma \ref{lemma-galvez-continuousfamily} assures the existence, for some \(\epsilon_0>0\), of a continuous family of radial solutions \(\{v_s\}_{s\in(c-\epsilon_0,c+\epsilon_0)}\) to the Weingarten equation on \(\overline{B_{\tilde R}}^C\) with zero boundary value, indexed by their leading coefficients \(s\)'s. We can assume that \(v_c\equiv v\). 
	
	For \(r\geq \tilde R\), let \(m_r\defeq\min_{\partial B_r}(u-v)\) and \(M_r\defeq\max_{\partial B_r}(u-v)\). Then, by the definition of leading coefficient and applying the maximum principle on sufficiently large annuli, we get:
	\[v_{c-\epsilon}+m_r\leq u\leq v_{c+\epsilon}+M_r\quad\text{ in }\overline{B_r}^C\]
	for every \(\epsilon\in(0,\epsilon_0)\). By continuity, we then find:
\begin{equation}\label{equation-harnack1}
		m_r\leq u-v\leq M_r\quad\text{ in }\overline{B_r}^C
\end{equation}
	for every \(r\geq \tilde R\). Now, we want to study how the difference \(M_r-m_r\) behaves as \(r\) varies, which will lead to the asymptotic \eqref{equation-asymptoticNEW}. We do this by means of the Harnack inequality for the linearized operator \eqref{equation-differentiated}. Consider the \(R_0>0\) given by Lemma \ref{lemma-differentiability}, and the asymptotic of \(A\) and \(b\) given b.  For \(r>2R_0\), observe that, on a ball of radius \(\frac{r}{2}\) centered on a point of \(\partial B_{r}\), the constant appearing in \cite[Corollary 9.25]{gilbargtrudinger} can be chosen, by what was said above, independent of \(r\). Now, pick \(R>2R_0\). By \eqref{equation-harnack1}, we know that \(w-m_{R}\) is nonnegative on \(B_{R}^C\). Moreover, since \(L\) is homogeneous, we have that \(L(w-m_{R})=0\). Then, we can apply Harnack inequality (see \cite[Corollary 9.25]{gilbargtrudinger}) to \(w-m_{R}\) on the circle \(\partial B_{2R}\), and get:
	\[M_{2R}-m_{R}\leq C(m_{2R}-m_{R}).\]
	Since \(m_r\) is non-decreasing in \(r\) and \(M_r\) is non-increasing in \(r\), we then get:
	\[M_{2R}-m_{2R}\leq\left(1-\frac{1}{C}\right)(M_{R}-m_{R})\]
	which implies, by iteration:
	\[M_{2^kR}-m_{2^kR}\leq\left(1-\frac{1}{C}\right)^k(M_{R}-m_{R})\]
	for any \(k\in\N\). Then, for $\alpha\defeq-\log_2(1-\frac{1}{C})$, we have:
	\[M_{2^kR}-m_{2^kR}\leq (2^k)^{-\alpha}(M_{R}-m_{R}),\]
	which implies, for any \(r\geq R\):
	\begin{equation*}
		M_r-m_r\leq \tilde Cr^{-\alpha}(M_{R}-m_{R}),
	\end{equation*}
	where \(\tilde C\) is a positive constant, only depending on \(R\) and, via \(C\), on $\Lambda$.
	
	Since \(m_r\) and \(M_r\) are both bounded and monotone, we then see that \(M_r,m_r\to b\), for some constant \(b\in\R\). Then:
	\begin{equation*}
		0\leq u(x)-v(x)-m_{|x|}\leq M_{|x|}-m_{|x|},
	\end{equation*}
	which readily implies \eqref{equation-asymptoticNEW}.
\end{proof}
\begin{remark}
	Observe that, in \eqref{equation-asymptoticNEW}, we actually have some control on the value of \(b\). Indeed, as in the proof of Theorem \ref{theorem-asymptotic}, restrict the domain, by choosing a bigger \(R\), so that the family of solutions with zero boundary value \(\{v_s\}_{s\in(c-\epsilon_0,c+\epsilon_0)}\) exists. Assume that \(v_c\equiv v\). Then, as in the proof of Theorem \ref{theorem-asymptotic}, we see that:
	\[\min_{\partial B_{R}}u\leq u-v\leq\max_{\partial B_R}u,\]
	that is, by \eqref{equation-asymptoticNEW}:
	\begin{equation*}
		\min_{\partial B_{R}}u\leq b\leq\max_{\partial B_R}u.
	\end{equation*}
\end{remark}
\begin{theorem}\label{theorem-uniqueness1}
	Let \(u\), \(v\) be unbounded solutions of \eqref{equation-weingarten}-\eqref{equation-uniformellipticity}-\eqref{equation-minimaltype} over $\overline{\Omega}^C$, where $\Omega$ is a smooth, bounded domain of \(\R^2\) such that $\overline{\Omega}^C$ is connected. Assume \(f'(0)\in(-1,0)\) that \(u\) and \(v\) have the same limit unit normal \(e_3\) and the same leading coefficient \(c>0\). Then:
	\begin{equation*}
		\min_{\partial \Omega}(u-v)\leq u-v\leq\max_{\partial\Omega}(u-v).
	\end{equation*}
	In particular, if \(u\) and \(v\) have the same boundary value, then \(u=v\).
\end{theorem}
\begin{proof}
	Let \(w\defeq u-v\). As we saw above, by Lemma \ref{lemma-differentiability}, there exists \(R_0>0\) such that \(w\) solves \eqref{equation-differentiated} on \(\overline{B_{R_0}}^C\). Define:
	\[M\defeq \max_{\partial B_{R_0}}w,\]
	and consider the function:
	\[\phi_\epsilon(x)\defeq M+\epsilon(|x|^\eta-R_0^\eta),\]
	for $\eta\in(0,1+f'(0))$.
	We first show that $\phi_\epsilon$ is a supersolution for \(L\). Indeed, by direct calculation, we find:
	\[L\left(|x|^\eta\right)=\eta |x|^{\eta-2}(\eta-1-f'(0)+o(1)),\]
	which is negative for \(|x|\) sufficiently big, by the choice of $\eta$. In particular, up to increasing \(R_0\) -- and changing the \(M\) accordingly -- we can assume that $\phi_\epsilon$ is a supersolution on \(\overline{B_{R_0}}^C\).
	
	Now observe that, by Theorem \ref{theorem-asymptotic}, we see that the difference \(u-v\) is bounded. On the other hand, we have that $\phi_\epsilon(x)\to+\infty$ as \(|x|\to+\infty\). Then, by comparing $\phi_\epsilon$ and \(w\) on sufficiently big annuli, standard maximum principle for \(L\) implies that:
	\[w\leq\phi_\epsilon\text{ on }\overline{B_{R_0}}^C.\]
	By the arbitrariness of $\epsilon>0$, we then see that:
	\[w\leq M=\max_{\partial B_{R_0}}w\text{ on }\overline{B_{R_0}}^C.\]
	The analogous argument applied to \(-w\) yields:
	\[w\geq\min_{\partial B_{R_0}}w\text{ on }\overline{B_{R_0}}^C.\]
	Then, the maximum of \(w\) over $\overline{\Omega}^C$ is attained in the compact set $\overline{B_{R_0}}\setminus\Omega$. More precisely, we can show that it is attained on $\partial\Omega$. Indeed, let \(L>0\) be sufficiently large that:
	\[u\leq v+L\quad\text{in }\overline{B_{R_0}}\setminus\Omega.\]
	Now, start lowering the graph of \(v|_{\overline{B_{R_0}}\setminus\Omega}\) until a first touching point with \(u|_{\overline{B_{R_0}}\setminus\Omega}\) occurs. If the touching point were not on $\partial\Omega$, then, by what was said above, we would have found an interior touching point between the graph of \(u\) and a vertical translate of the graph of \(v\), with the normal vectors to the graph pointing in the same direction. But this would lead to a contradiction by the maximum principle for the elliptic Weingarten equation. An analogous argument shows that the minimum of \(w\) is also attained on\,\,$\partial\Omega$.
\end{proof}
The following maximum principle at infinity follows:
\begin{corollary}\label{corollary-maximumprincipleatinfinity}
	Let \(u\), \(v\) be unbounded solutions of \eqref{equation-weingarten}-\eqref{equation-uniformellipticity}-\eqref{equation-minimaltype} over the exterior of some \(B_R\), with \(u,v>0\) and having the same limit unit normal \(e_3\). Assume \(f'(0)\in(-1,0)\), \(u\geq v\) in \(B_R^C\) and that there exists a diverging sequence \(\{x_n\}_{n\in\N}\subset B_R^C\) such that:
	\[\lim_{n\to \infty}u(x_n)-v(x_n)=0.\]
	Then \(u\equiv v\).
\end{corollary}
\begin{proof}
	From the existence of the sequence \(\{x_n\}_{n\in\N}\), we have that \(u\) and \(v\) have the same leading coefficient \(c\). Moreover, applying the maximum principle in the bounded annulus \(B_{3R}\setminus \overline{B_R}\), we have that, if \(u\) and \(v\) do not coincide, then \(u-v>0\) on \(\partial B_{2R}\). That is, by compactness:
	\[\min_{\partial B_{2R}}(u-v)\geq\delta>0.\]
	Applying Theorem \ref{theorem-uniqueness1} in \(\overline{B_{2R}}^C\), this contradicts the existence of the sequence \(\{x_n\}_{n\in\N}\). 
\end{proof}
\begin{remark}\label{remark-maximumprincipleatinfinity}
	Observe that the maximum principle at infinity cannot hold for ends of finite limiting height with \(f'(0)\neq-1\). Indeed, by Lemma \ref{lemma-galvez-continuousfamily}, one can first pick a solution \(v\) satisfying:
	\[\begin{cases}
		\lim_{x\to+\infty}v(|x|)=h&\text{ in }\overline{B_1}^C,
		\\
		v|_{\partial B_1}=0,
	\end{cases}\]
	for some \(h>0\). Now, consider the constant function \(u\equiv h\) in \(\overline{B_1}^C\). By the standard maximum principle we have \(u\geq v\) in \(B_1^C\), but \(u(x)-v(x)\to0\) as \(|x|\to+\infty\).
\end{remark}
Even though the maximum principle at infinity, as stated in Corollary \ref{corollary-maximumprincipleatinfinity}, cannot hold for bounded ends, the following comparison principle holds for every class of ends:
\begin{theorem}\label{theorem-weakmaximumprincipleatinfinity}
	Let \(u\), \(v\) be solutions of \eqref{equation-weingarten}-\eqref{equation-uniformellipticity}-\eqref{equation-minimaltype} over the exterior of some bounded domain \(\Omega\) such that $\overline{\Omega}^C$ is connected, with \(u,v\) bounded from below and having the same limit unit normal \(e_3\) and same boundary value. Assume that one of the following holds:
	\begin{itemize}
		\item \(u\) and \(v\) have the same limiting height \(h\) and, in the case \(h=+\infty\), they have the same leading coefficient \(c\);
		\item There is a diverging sequence \(\{x_n\}_{n\in\N}\subset\Omega^C\) such that:
		\[\lim_{n\to \infty}u(x_n)-v(x_n)=0.\]
	\end{itemize}
	Then \(u\equiv v\).
\end{theorem}
\begin{proof}
	It is easy to see that the existence of the diverging sequence \(\{x_n\}_{n\in\N}\) implies that \(u\), \(v\) have the same limiting height and, in the case this is infinite, the same leading coefficient. Let us then assume that \(u\) and \(v\) have the same limiting height and the same leading coefficient.
	
	The proof is divided into two cases:
	\\
	\fbox{\textbf{Bounded ends}}\,: In this case, we have that \(u\) and \(v\) are solutions to $\mathcal{P}(\overline{\Omega}^C,\phi)$ with the same limiting height \(h\in\R\). Then, if \(u\) and \(v\) did not coincide, we could find a vertical translate of the graph of \(v\) touching the graph of \(u\) at an interior point, thereby leading to a contradiction by the maximum principle.
	\\
	\fbox{\textbf{Unbounded ends}}\,: In the case \(f'(0)\neq-1\), it is an immediate consequence of Theorem \ref{theorem-uniqueness1}. In any case, we can prove it using maximum principle at infinity -- that is, Corollary \ref{corollary-maximumprincipleatinfinity} for \(f'(0)\neq-1\) and \cite[Corollary 1.5]{galvezbarbieri2026asymptotic} for \(f'(0)=-1\). Indeed, by the asymptotic expansions \eqref{equation-asymptoticNEW} (for \(f'(0)\neq-1\)) and \eqref{e1.10} (for \(f'(0)=-1\)) there exists a vertical translate of \(v\), say \(v+L\), for some \(L>0\), which is entirely above \(u\):
	\[
	v+L>u\text{ in }\Omega^C.
	\]
	We can start lowering the graph of \(v+L\), by decreasing \(L\), until either we find a first contact point between the graph of \(u\) and the graph of \(v+L\), or a sequence \(\{x_n\}_n\subset\Omega^C\) with \(|x_n|\to+\infty\) and:
	\[v(x_n)+L-u(x_n)\to0.\]
	In any case, if this happened for a positive \(L\), standard maximum principle or maximum principle at infinity would lead to a contradiction, since the two graphs should then coincide, with different boundary values. Then, we find that the first touching point between the two graphs must happen at the boundary, and for \(L=0\). This implies that \(v\geq u\).
	
	On the other hand, if we start by translating \(u\), we get \(u\geq v\). Then \(u\) must coincide with \(v\).
\end{proof}
As a corollary, we get the following symmetry inheritance result. An analogous property for Weingarten surfaces of cmc type has been investigated in \cite{angelo}, where we studied properly embedded Weingarten surfaces (not necessarily graphs) of cmc type with boundary on a strictly convex planar curve.
\begin{corollary}\label{corollary-symmetry}
	Let $\Omega\subset\R^2$ be a bounded domain such that $\overline{\Omega}^C$ is connected, and let \(S\) be an isometry of \(\R^2\) with \(S(\Omega)=\Omega\). Consider a continuous function \(\phi:\partial\Omega\to\R\) such that \(\phi\circ S\equiv\phi\) and assume that there exists a solution \(u\) to $\mathcal{P}(\overline{\Omega}^C,\phi)$, having limit unit normal \(e_3\). Then \(u\circ S\equiv u\).
\end{corollary}
\begin{proof}
	First of all observe that, by Remark \ref{remark-unboundedends}, we can assume \(u\) to be bounded from below. Now, observe that the map:
	\[(x,y,z)\to(S(x,y),z)\]
	is an isometry of \(\R^3\) and therefore it preserves the principal curvatures of the graph of \(u\). Then, \(u\circ S\) is still a solution to $\mathcal{P}(\overline{\Omega}^C,\phi)$ with the same limiting height and, when this is infinite, same leading coefficient. Therefore Theorem \ref{theorem-weakmaximumprincipleatinfinity} applies.
\end{proof}
\section{The foliation property}\label{section-krust}
Once the existence of a solution to the Dirichlet problem is proved, one may want to address uniqueness properties. But, for the exterior Dirichlet problem, one cannot expect uniqueness of the solution in general. Indeed, already for $\mathcal{P}(\overline{B_1}^C,0)$, uniqueness fails dramatically: rotational solutions (which, by Corollary \ref{corollary-symmetry}, are the only solutions with limit unit normal \(e_3\) to the problem) form a one-parameter family. This is what Lemma \ref{lemma-galvez-continuousfamily} basically says.

In this section, we prove that this property of having a one-parameter family of solutions is, in a sense, not merely a feature of $\mathcal{P}(\overline{B_1}^C,0)$. Let \(\Omega\) be a bounded domain with \(C^0\) boundary, such that $\overline{\Omega}^C$ is connected, and let \(\phi\in C^0(\partial\Omega)\). Assume that there exist two distinct solutions \(u_a\), \(u_b\) to $\mathcal{P}(\overline{\Omega}^C,\phi)$, bounded from below and with the same limit unit normal \(e_3\). Then, we will show that the set:
\[\left\{\text{graph}(u)\,:\,u\text{ is a solution to }\mathcal{P}(\overline{\Omega}^C,\phi),\,\,u_a\leq u\leq u_b\right\}\]
forms a foliation of the region between the graphs of \(u_a\) and \(u_b\).
\begin{theorem}\label{theorem-foliation1}
	Let \(u_a\), \(u_b\) be distinct solutions to $\mathcal{P}(\overline{\Omega}^C,\phi)$, bounded from below and with limit unit normal \(e_3\). Let \(h_a\), \(h_b\), \(c_a\), \(c_b\) be the limiting heights and leading coefficients of \(u_a\) and \(u_b\) respectively, and assume that \(u_a\leq u_b\). Then, the set of solutions to $\mathcal{P}(\overline{\Omega}^C,\phi)$ bounded between \(u_a\) and \(u_b\) is a foliation of the region: \[\left\{(x_1,x_2,x_3)\,:\,u_a(x_1,x_2)<x_3<u_b(x_1,x_2)\right\}.\]
	These are given by the following (possibly empty) two one-parameter families of solutions:
	\begin{itemize}
		\item \(\{u_c\}_{c\in(c_a,c_b)}\), a family of unbounded solutions, indexed by their leading coefficients. This family is empty if and only if \(h_b<+\infty\);
		\item \(\{w_h\}_{h\in (h_a,h_+]}\), for some constant \(h_+\), a family of bounded solutions, indexed by their limiting heights, that is:
		\[h\defeq\lim_{|x|\to+\infty}w_h(x).\]
		This family is empty if and only if the first family is not empty and \(u_a\) can be obtained as the limit of the\,\,\(u_c\)'s as \(c\to c_a^+\).
	\end{itemize} 
\end{theorem}
\begin{proof}
	We first study the two mutually exclusive cases \(h_a=+\infty\) and \(h_b<+\infty\). In the first case, we construct the family \(\{u_c\}_{c\in(c_a,c_b)}\), and in the second case we construct the family \(\{w_h\}_{h\in(h_a,h_b)}\). The general case will be solved in the end.
	
	The idea is to construct, for each possible value of \(c\in(c_a,c_b)\) (resp. \(l\in(h_a,h_b)\)), a solution with leading coefficient \(c\) (resp. limiting height \(l\)). The construction is similar to the one of Section \ref{section-exteriordirichlet}: we solve the Dirichlet problem on bounded regions exhausting \(\overline{\Omega}^C\) using the Perron method, and then use the convergence results of Section \ref{section-dirichlet} to construct a global solution, with the specified asymptotic. Observe that, since every solution that we are going to build will be below \(u_b\), then, by Theorem \ref{theorem-gradientestimate}, every leaf of the foliation will have limit unit normal \(e_3\).
	
	Let \(n\in\N\) sufficiently big, so that $\overline{\Omega}\subset B_n(0)$, and define the annular region:
	\[A_n\defeq B_n(0)\setminus\overline{\Omega}.\]
	
	\fbox{\textbf{Case 1: \(h_a=+\infty\)}}\,: From Theorem \ref{theorem-weakmaximumprincipleatinfinity}, we know that \(0<c_a<c_b\). Fix \(c\in(c_a, c_b)\). We want to solve a Dirichlet problem on \(A_n\). The boundary value on $\partial\Omega$ will be equal to \(u_a|_{\partial\Omega}=u_b|_{\partial\Omega}\). Then, \(u_a\) and \(u_b\) constitute a barrier at every point of $\partial\Omega$. Regarding the boundary value on $\partial B_n$, we choose it as follows. By Lemma \ref{lemma-galvez-continuousfamily}, there exists a radius \(R\) and a radial solution \(w\) to \(\mathcal{P}(\overline{B_R(0)}^C,0)\) with leading coefficient \(c\). For \(n>R\), the boundary value on $\partial B_n(0)$ will be the constant \(w|_{\partial B_n(0)}\). As in Section \ref{section-exteriordirichlet}, at every point of \(\partial B_n(0)\), for \(n\) sufficiently big, an horizontal plane and a tilted plane define a barrier. Then, Proposition \ref{proposition-perron} can be applied: there exists \(v_n\in C^2(A_n)\cap C^0(\overline{A_n})\), solving the Weingarten equation, such that:
	\[v_n|_{\partial\Omega}\equiv u_a|_{\partial\Omega}(\equiv u_b|_{\partial\Omega})\,\,\text{ and }\,\,v_n|_{\partial B_n(0)}\equiv w|_{\partial B_n(0)}.\]
	Now, the results of Section \ref{section-dirichlet} can be used to extract from \(\{v_n\}_{n\in\N}\) a subsequence converging in \(C^2_{\text{loc}}(\overline{\Omega}^C)\) to a solution \(u_c\) to the Weingarten equation in $\overline{\Omega}^C$. Moreover, by construction, we have that \(u_a\leq u_c\leq u_b\) in $\overline{\Omega}^C$, and this implies that \(u_c\in C^2(\overline{\Omega}^C)\cap C^0(\Omega^C)\), with:
	\[u_c|_{\partial\Omega}\equiv u_a|_{\partial\Omega}.\]
	Regarding the leading coefficient of \(u_c\), observe that \(u_c\), outside of a compact set, is trapped between two vertical translates of \(w\), which implies that the leading coefficient of \(u_c\) is in fact \(c\). Moreover, by Theorem \ref{theorem-weakmaximumprincipleatinfinity}, \(u_c\) is the only solution to \(\mathcal{P}(\overline{\Omega}^C,\phi)\) with leading coefficient \(c\), and therefore \(u_c\) is well defined.

	To finish the proof, we need to prove that the set of the graphs of the \(u_c\)'s form in fact a foliation of the region:
	\[\left\{(x_1,x_2,x_3)\,:\,(x_1,x_2)\in\Omega^C,\,\,u_a(x_1,x_2)<x_3<u_b(x_1,x_2)\right\}.\]
	That is, we need to prove that, for every \((x_1,x_2,x_3)\) such that \(u_a(x_1,x_2)<x_3<u_b(x_1,x_2)\), there exists a unique \(c\in(c_a,c_b)\) such that:
	\[u_c(x_1,x_2)=x_3.\]
	
	For \((x_1,x_2,x_3)\) in the set, define:
	\[c\defeq\sup\left\{s\in(c_a,c_b)\,:\,u_s(x_1,x_2)<x_3\right\},\]
	and observe that, necessarily, one has:
	\[c=\inf\left\{s\in(c_a,c_b)\,:\,u_s(x_1,x_2)>x_3\right\}.\]
	Then, by the results of Section \ref{section-dirichlet}, \(u_c\) can be constructed as the \(C^2_\text{loc}(\overline{\Omega}^C)\) limit of two sequences of functions: one of them gives \(u_c(x_1,x_2)\leq x_3\) and the other one gives \(u_c(x_1,x_2)\geq x_3\). Then \(u_c(x_1,x_2)=x_3\). The uniqueness of the leaf comes from Theorem \ref{theorem-weakmaximumprincipleatinfinity}.
	
	\fbox{\textbf{Case 2: \(h_b<+\infty\)}}\,: Boundedness of \(u_a\) and \(u_b\) requires this case to be treated separately. We still start by solving the Dirichlet problem on \(A_n\). Let \(h\in(h_a, h_b)\). The boundary datum on $\partial\Omega$ is still \(u_a|_{\partial\Omega}\), and \(u_a\), \(u_b\) are barriers at every point of $\partial\Omega$. Regarding $\partial B_n(0)$, observe that we want to construct a solution that converges to the horizontal plane at height \(h\). Then, the boundary datum that we assign on $\partial B_n(0)$ is simply the constant \(h\). As barriers at every point of \(\partial B_n(0)\), we pick two suitable tilted planes -- observe that, unlike the preceding case in which one of the two planes could be chosen horizontal, at least for big \(n\), in this case the boundedness of the solutions implies that a horizontal plane could fail to be a supersolution!
	
	We then find a suitable solution on each \(A_n\) and, by compactness, we get a function \(w_h\in C^2(\overline{\Omega}^C)\cap C^0(\Omega^C)\) solving the Weingarten equation, with boundary datum \(u_a|_{\partial\Omega}\). How to prove that \(w_h\) has the right limiting height? In this case, we cannot use radial solutions to trap \(w_h\), since a vertical translation also changes the limiting height of the solution. This same issue can be exploited as follows. Since \(h_a<h<h_b\), we can trap \(w_h\) between \(u_b-(h_b-h)-\epsilon\) and \(u_a+(h-h_a)+\epsilon\), for $\epsilon>0$. Then, by standard maximum principle applied to the solutions of the local Dirichlet problems, we find:
	\[u_b-(h_b-h)-\epsilon\leq w_h\leq u_a+(h-h_a)+\epsilon\]
	and, by arbitrariness of $\epsilon>0$, we get:
	\[\lim_{|x|\to+\infty}u_h(x)=h.\]
	That is, \(h\) is the limiting height of \(w_h\). To prove that the graphs of the \(w_h\)'s form a foliation, we can argue as above.
	
	\fbox{\textbf{General case}}\,: We only need to prove the theorem when \(u_a\) has finite limiting height and \(u_b\) has infinite limiting height:
	\[0=c_a<c_b.\]
	In this case, the same proof as in Case 1 gives the existence of the family \(\{u_c\}_{c\in(0,c_b)}\).
	\begin{figure}[!htb]
		\centering
		\includegraphics[trim={4cm 18cm 2cm 4cm}, clip, width=\textwidth]{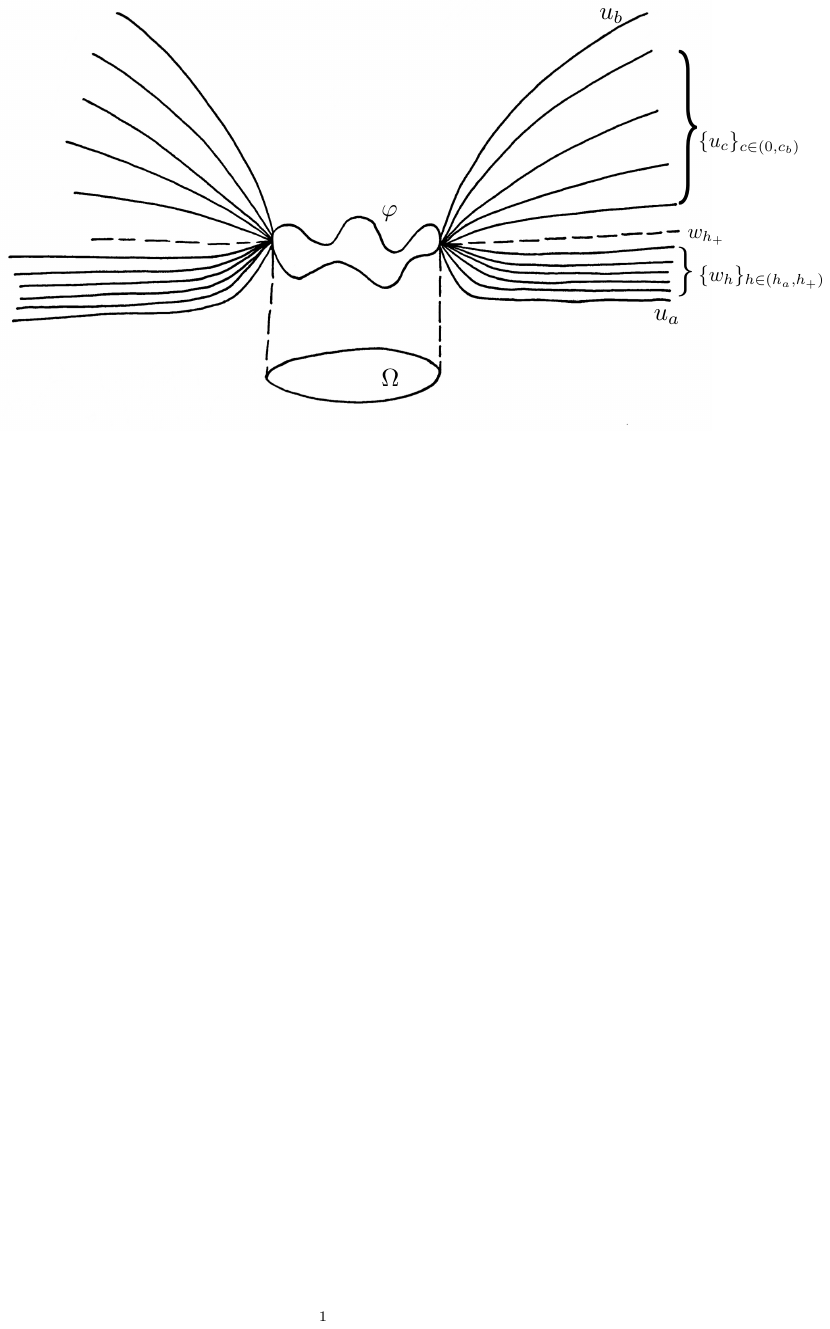}
		\caption{The case in which \(0=c_a<c_b\)}
	\end{figure}
	Now, by the results of Section \ref{section-dirichlet}, a subsequence of the family converges, as \(c\to0^+\), to a function \(u\) solving $\mathcal{P}(\overline{\Omega}^C,\phi)$, whose leading coefficient is easily seen to be 0. If \(u\equiv u_a\), there is nothing else to prove. On the other hand, if this is not the case, then Theorem \ref{theorem-weakmaximumprincipleatinfinity} implies that:
	\[\lim_{|x|\to+\infty}u_a\eqdef h_-<h_+\defeq\lim_{|x|\to+\infty}u\]
	and we argue as above to construct the family \(\{w_h\}_{h\in(h_-,h_+]}\). The foliation property is then proved separately in the two regions:
	\[\left\{(x_1,x_2,x_3)\,:\,u(x_1,x_2)<x_3<u_b(x_1,x_2)\right\}\]
	and:
	\[\left\{(x_1,x_2,x_3)\,:\,u_a(x_1,x_2)<x_3<u(x_1,x_2)\right\}\]
	as in the preceding cases.
\end{proof}
\textbf{Acknowledgments:} This study was partially supported by INdAM-GNSAGA. The author is grateful to Barbara Nelli for first proposing the problem and for many fruitful discussions, and to José Antonio Gálvez for stimulating his interest in the theory of elliptic Weingarten surfaces. The paper was written while the author was visiting IMPA, whose warm hospitality is sincerely appreciated. Part of the stay at IMPA was partially supported by INdAM and IMPA, under the INdAM-IMPA Call 2026 for scientific visits in Brazil.

\textbf{AI disclosure:} In Section \ref{section-maximumprincipleinfinity}, the strategy to differentiate the Weingarten equation and the right supersolution emerged from interactions between the author and a LLM. All the mathematical content was independently verified and reworked by the author. The work does not contain any AI written text.

\end{document}